\documentclass{article}
\usepackage{graphicx}
\usepackage[utf8]{inputenc}
\usepackage[english]{babel}
\usepackage[authoryear,sort]{natbib}
\usepackage{fullpage}
\usepackage{color}
\usepackage{enumerate}
\usepackage{multicol}
\usepackage{amsthm,amsmath,amssymb}
\usepackage[linesnumbered,ruled,vlined]{algorithm2e}
\newtheorem{lemma}{Lemma}
\newtheorem{thm}{Theorem}
\newcounter{prf}
\theoremstyle{definition}
\newtheorem{definition}{Definition}

\title{Selecting optimal minimum spanning trees that share a topological correspondence with phylogenetic trees.}
\author{Prabhav Kalaghatgi\\
\small{Max Planck Institute for Informatics}\\
\small{Saarbr\"{u}cken}\\
prabhavk@mpi-inf.mpg.de
\and
Thomas Lengauer\\
\small{Max Planck Institute for Informatics}\\
\small{Saarbr\"{u}cken}\\
lengauer@mpi-inf.mpg.de
}
\date {}
\begin{document}
\maketitle
\begin{abstract}
\cite{Choi2010b} introduced a minimum spanning tree (MST)-based method called CLGrouping, for constructing tree-structured probabilistic graphical models, a statistical framework that is commonly used for inferring phylogenetic trees. While CLGrouping works correctly if there is a unique MST, we observe an indeterminacy in the method in the case that there are multiple MSTs. In this work we remove this indeterminacy by introducing so-called vertex-ranked MSTs. We note that the effectiveness of CLGrouping is inversely related to the number of leaves in the MST. This motivates the problem of finding a vertex-ranked MST with the minimum number of leaves (MLVRMST). We provide a polynomial time algorithm for the MLVRMST problem, and prove its correctness for graphs whose edges are weighted with tree-additive distances.
\end{abstract}
\section{Introduction}
Phylogenetic trees are commonly modeled as tree-structured probabilistic graphical models with two types of vertices: labeled vertices that represent observed taxa, and hidden vertices that represent unobserved ancestors. The length of each edge in a phylogenetic tree quantifies evolutionary distance. If the set of taxa under consideration contain ancestor-descendant pairs, then the phylogenetic tree has labeled internal vertices, and is called a \emph{generally labeled tree} \citep{Kalaghatgi2016}. The data that is used to infer the topology and edge lengths is usually available in the form of gene or protein sequences.

Popular distance-based methods like neighbor joining (NJ; \cite{Saitou1987}) and BIONJ \citep{Gascuel1997} construct phylogenetic trees from estimates of the evolutionary distance between each pair of taxa. \cite{Choi2010b} introduced a distance-based method called Chow-Liu grouping (CLGrouping). \cite{Choi2010b} argue that CLGrouping is more accurate than NJ at reconstructing phylogenetic trees with large diameter. The diameter of tree is the number of edges in the longest path of the tree.

CLGrouping operates in two phases. The first phase constructs a distance graph $G$ which is a complete graph over the labeled vertices where each edge is weighted with the distance between each pair of labeled vertices. Subsequently a minimum spanning tree (MST) of $G$ is constructed. In the second phase, for each internal vertex $v_{i}$ of the MST, the vertex set $V_{i}$ consisting of $v_{i}$ and its neighbors is constructed. Subsequently a generally labeled tree $T_{i}$ over $V_{i}$ is inferred using a distance-based tree construction method like NJ. The subtree in the MST that is induced by $V_{i}$ is replaced with $T_{i}$.

Distances are said to be additive in a tree $T$ if the distance between each pair of vertices $u$ and $v$ is equal to the sum of lengths of edges that lie on the path in $T$ between $u$ and $v$. Consider the set of all phylogenetic trees $\mathcal{T}$ such that the edge length of each edge in each tree in $\mathcal{T}$ is strictly greater than zero. A distance-based tree reconstruction method is said to be consistent if for each $\{D,T|T\in \mathcal{T}\}$ such that $D$ is additive in $T$, the tree that is reconstructed using $D$ is identical to $T$. Please note the following well-known result regarding the correspondence between trees and additive distances. Considering all trees in $\mathcal{T}$, if $D$ is additive in a tree $T$ then $T$ is unique \citep{Buneman1971}. 

We show that if $G$ has multiple MSTs then CLGrouping is not necessarily consistent. We show that there always exists an MST $M$ such that CLGrouping returns the correct tree when $M$ is used in the second phase of CLGrouping. We show that $M$ can be constructed by assigning ranks to the vertices in $G$, and by modifying standard MST construction algorithms such that edges are compared on the basis of both edge weight and ranks of the incident vertices. The MSTs that are constructed in this manner are called vertex-ranked MSTs.

Given a distance graph, there may be multiple vertex-ranked MSTs with vastly different number of leaves. \cite{Huang2014} showed that CLGrouping affords a high degree of parallelism, because, phylogenetic tree reconstruction for each vertex group can be performed independently. With respect to parallelism, we define an optimal vertex-ranked MST for CLGrouping to be a vertex-ranked MST with the maximum number of vertex groups, and equivalently, the minimum number of leaves.

We developed an $O(n^{2}\log n)$ time algorithm Algo. \ref{algo:MLVRMSTFromDistanceGraph} that takes as input a distance graph and outputs a vertex-ranked MST with the minimum number of leaves (MLVRMST). The proof of correctness of Algo. \ref{algo:MLVRMSTFromDistanceGraph} assumes that the edges in the distance graph are weighted with tree-additive distances. 

\section{Terminology}
A phylogenetic tree is an undirected edge-weighted acyclic graph with two types of vertices: labeled vertices that represent observed taxa, and hidden vertices that represent unobserved taxa. Information, e.g., in the form of genomic sequences, is only present at labeled vertices. We refer to the edge weights of a phylogenetic tree as edge lengths. The length of an edge quantifies the estimated evolutionary distance between the sequences corresponding to the respective incident vertices. All edge lengths are strictly positive. Trees are leaf-labeled if all the labeled vertices are leaves. Leaf-labeled phylogenetic trees are the most commonly used models of evolutionary relationships. Generally labeled trees are phylogenetic trees whose internal vertices may be labeled, and are appropriate when ancestor-descendant relationships may be present in the sampled taxa \citep{Kalaghatgi2016}.

Each edge in a phylogenetic tree partitions the set of all labeled vertices into two disjoint sets which are referred to as the split of the edge. The two disjoint sets are called to the sides of the split. 

A phylogenetic tree can be rooted by adding a hidden vertex (the root) to the tree, removing an edge $e$ in the tree, and adding edges between the root and the vertices that were previously incident to $e$. Edge lengths for the newly added edges must be positive numbers and must sum up to the edge length of the previously removed edge. Rooting a tree constructs a directed acyclic graph in which each edge is directed away from the root. 

A leaf-labeled phylogenetic tree is clock-like if the tree can be rooted in such a way that all leaves are equidistant from the root. Among all leaf-labeled phylogenetic trees, maximally balanced trees and caterpillar trees have the smallest and largest diameter, respectively, where the diameter of a tree is defined as the number of edges along the longest path in the tree. 

The distance graph $G$ of a phylogenetic tree $T$ is the edge-weighted complete graph whose vertices are the labeled vertices of $T$. The weight of each edge in $G$ is equal to the length of the path in $T$ that connects the corresponding vertices that are incident to the edge. A minimum spanning tree (MST) of an edge-weighted graph is a tree that spans all the vertices of the graph, and has the minimum sum of edge weights.
\section{Chow-Liu grouping}
\citet{Choi2010b} introduced the procedure Chow-Liu grouping (CLGrouping) for the efficient reconstruction of phylogenetic trees from estimates of evolutionary distances. If the input distances are additive in the phylogenetic tree $T$ then the authors claim that CLGrouping correctly reconstructs $T$. 

CLGrouping consists of two stages. In the first stage, an MST $M$ of $G$ is constructed. In the second stage, for each internal vertex $v$, a vertex group $Nb(v)$ is defined as follows: $Nb(v)$ is the set containing $v$ and all the vertices in $M$ that are adjacent to $v$. For each vertex group, a phylogenetic tree $T_{v}$ is constructed using distances between vertices in $Nb(v)$. Subsequently, the graph in $M$ that is induced by $Nb(v)$ is replaced by $T_{v}$ (see Fig. \ref{fig:CLGroupingUsingSortedAndUnsortedDistances}e for an illustration). $T_{v}$ may contain hidden vertices which may now be in the neighborhood of an internal vertex $w$ that has not been visited as yet. If this the case, then we need an estimate of the distance between the newly introduced hidden vertices and vertices in $Nb(w)$. Let $h_{v}$ be the hidden vertex that was introduced when processing the internal vertex $v$. The distance from $h_{v}$ to a vertex $k \in Nb(w)$ is estimated using the following formula, $d_{h_{v}k} = d_{vk} - d_{vh_{v}}$.

The order in which the internal vertices are visited is not specified by the authors and does not seem to be important. CLGrouping terminates once all the internal vertices of $M$ have been visited.

This procedure is called Chow-Liu grouping because the MSTs that are constructed using additive distances are equivalent to Chow-Liu trees \citep{Chow1968}, for certain probability distributions. Please read \cite{Choi2010b} for further detail.


\section{Indeterminacy of CLGrouping}

\begin{figure}[!t]
\centering
\includegraphics[width=0.6\textwidth]{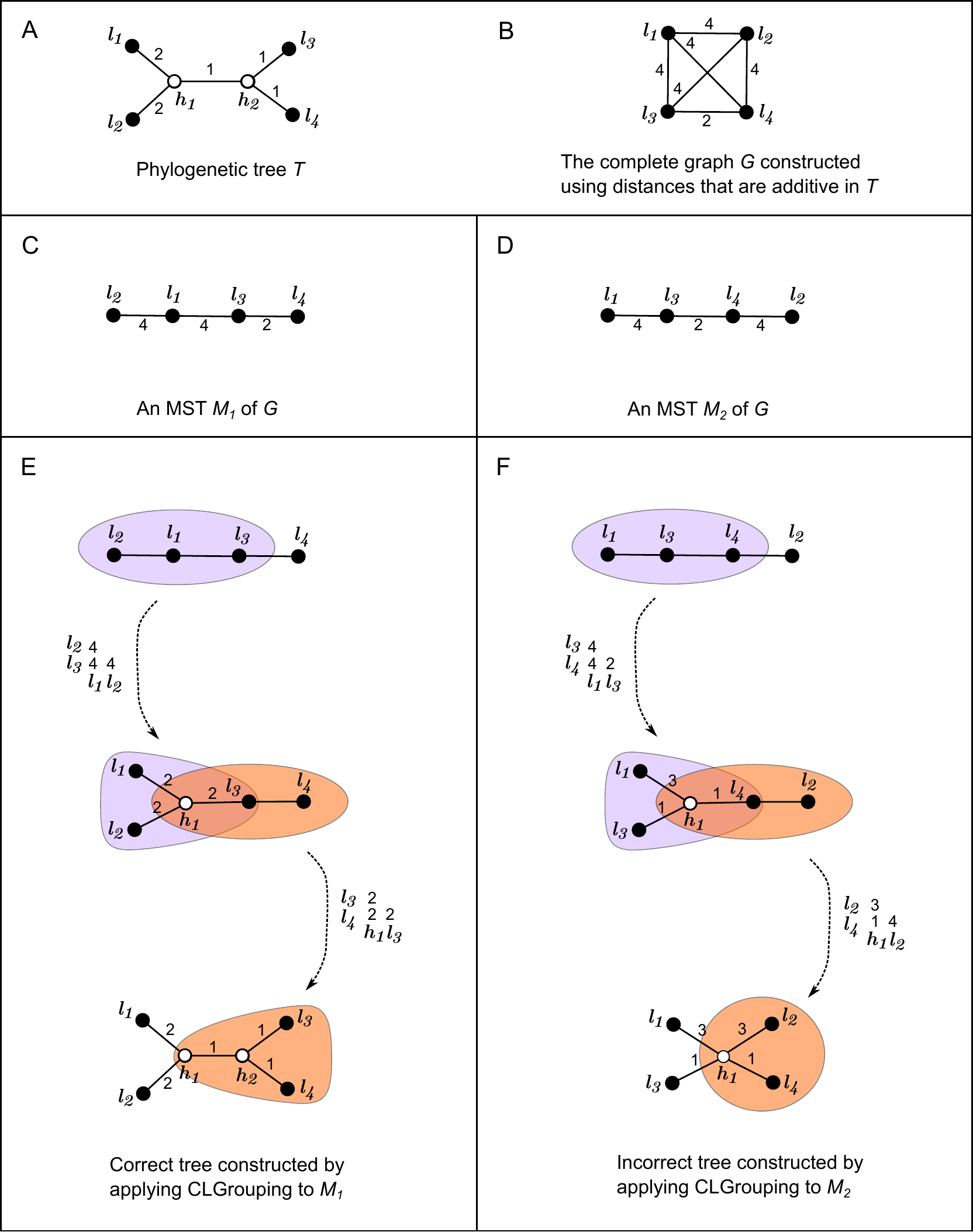}
\caption{The example used to demonstrate that CLGrouping may not reconstruct the correct tree if there are multiple MSTs. The phylogenetic tree $T$ that is used in this example is shown in panel a. The distance graph $G$ of $T$ is shown in panel b. Two MSTs of $G$, $M_{1}$ and $M_{2}$, respectively, are shown in panels c and d. Panels e and f show the intermediate steps and the final result	of applying CLGrouping to $M_{1}$ and $M_{2}$ respectively. CLGrouping reconstructs the original phylogenetic tree if it is applied to $M_{1}$, but not if it is applied to $M_{2}$.}
\label{fig:CLGroupingUsingSortedAndUnsortedDistances}
\end{figure}

CLGrouping is not necessarily consistent if there are multiple MSTs. We demonstrate this with the phylogenetic tree $T$ shown in Fig. \ref{fig:CLGroupingUsingSortedAndUnsortedDistances}a. For the corresponding distance graph $G$ of $T$ (see Fig. \ref{fig:CLGroupingUsingSortedAndUnsortedDistances}b), two MSTs of $G$, $M_{1}$ and $M_{2}$ are shown in Fig. \ref{fig:CLGroupingUsingSortedAndUnsortedDistances}c and Fig. \ref{fig:CLGroupingUsingSortedAndUnsortedDistances}d, respectively. The intermediate steps, and the final result of applying CLGrouping to $M_{1}$ and $M_{2}$ are shown in Fig. \ref{fig:CLGroupingUsingSortedAndUnsortedDistances}e and Fig. \ref{fig:CLGroupingUsingSortedAndUnsortedDistances}f, respectively. CLGrouping reconstructs the original phylogenetic tree if it is applied to $M_{1}$ but not if it is applied to $M_{2}$.

The notion of a surrogate vertex is central to proving the correctness of CLGrouping. The surrogate vertex of a hidden vertex is the closest labeled vertex, w.r.t. distances defined on the phylogenetic tree. CLGrouping will reconstruct the correct phylogenetic tree only if the MST can be constructed by contracting all the edges along the path between each hidden vertex and its surrogate vertex. Since the procedure that constructs the MST is not aware of the true phylogenetic tree, the surrogate vertex of each hidden vertex must selected implicitly. In the example shown earlier, $M_{1}$ can be constructed by contracting the edges $(h_{1},l_{1})$, and $(h_{2},l_{3})$. Clearly there is no selection of surrogate vertices such that $M_{2}$ can be constructed by contracting the path between each hidden vertex and the corresponding surrogate vertex.

If there are multiple labeled vertices each of which is closest to a hidden vertex then \cite{Choi2010b} assume that the corresponding surrogate vertex is implicitly selected using the following tie-breaking rule.

\begin{figure}[!ht]
\centering
\includegraphics[width=0.2\textwidth]{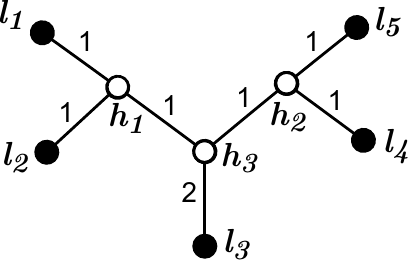}
\caption{The phylogenetic tree that is used to demonstrate that the tie-breaking rule as defined by \citet{Choi2010b} cannot be applied in general.}
\label{fig:TieBreakingRuleCannotBeGenerallyApplied}
\end{figure}

Let the surrogate vertex set $\mathbf{Sg}(h)$ of a vertex $h$ be the set of all labeled vertices that are closest to $h$. If $l_{1}$ and $l_{2}$ belong to both $\mathbf{Sg}({h_{1}})$ and $\mathbf{Sg}({h_{1}})$, then the same labeled vertex (either $l_{1}$ or $l_{2}$) is selected as the surrogate vertex of both $h_{1}$ and $h_{2}$. This rule for selecting surrogate vertices cannot be consistently applied across all hidden vertices. We demonstrate this with an example.  For the tree shown in Fig. \ref{fig:TieBreakingRuleCannotBeGenerallyApplied} we have $\mathbf{Sg}(h_{1})=\{l_{1},l_{2}\}$, $\mathbf{Sg}(h_{2})=\{l_{4},l_{5}\}$, and $\mathbf{Sg}(h_{3}$)=$\{l_{1},l_{2},l_{3},l_{4},l_{5}\}$. It is clear that there is no selection of surrogate vertices that satisfies the tie-breaking rule.

\section{Ensuring the consistency of CLGrouping}

In order to construct an MST that is guaranteed to have the desired topological correspondence with the phylogenetic tree, we propose the following tie-breaking rule for selecting the surrogate vertex. Let there be a total order over the set of all labeled vertices. Let $\mathcal{R}(l)$ be the rank of vertex $l$ that is given by the order. We define the surrogate vertex $\mbox{Sg}(h)$ of $h$ to be the highest ranked labeled vertex among the set of labeled vertices that are closest to $h$. That is,
\begin{definition}
\label{def:vertexRankBasedSurrogateVertex}
$$\mbox{Sg}(h) = \min\limits_{l\in\mathbf{Sg}(h)}\mathcal{R}(l)\mbox{ ,where,}$$
$$\mathbf{Sg}(h) = \min_{l \in \mathcal{L}(T)}d_{lh}.$$
\end{definition}

The inverse surrogate set $\mbox{Sg}^{-1}(l)$ is the set of all hidden vertices whose surrogate vertex is $l$.

In order to ensure that the surrogate vertices are selected on the basis of both distance from the corresponding hidden vertex and vertex rank, it is necessary that information pertaining to vertex rank is used when selecting the edges of the MST. We use Kruskal's algorithm \citep{Kruskal1956} for constructing the desired MST. Since Kruskal's algorithm takes as input a set of edges sorted w.r.t. edge weight, we modify the input by sorting edges with respect to edge weight and vertex rank as follows. It is easy to modify other algorithms for constructing MSTs in such a way that vertex rank is taken into account.

\begin{definition}
\label{def:orderedEdges}
We define below, what is meant by sorting edges on the basis of edge weight and vertex rank. Given a edge set $E$, and a ranking $\mathcal{R}$ over vertices in $E$, let $d(u,v)$ be the weight of the edge $\{u,v\}$, and let $\mathcal{R}(u)$ be the rank of the vertex $u$. Let the relative position of each pair of edges in the list of sorted edges be defined using the total order $<$. That is to say, for each pair of edges $\{a,b\}$ and $\{c,d\}$,
\begin{align*}
&\{a,b\} < \{c,d\}\mbox{, if and only if}\\
&(i)\,\,\,\,\,d(a,b) < d(c,d),\mbox{ or if}\\
&(ii)\,\,\,d(a,b) = d(c,d)\mbox{ and }\min\{\mathcal{R}(a),\mathcal{R}(b)\} < \min\{\mathcal{R}(c),\mathcal{R}(d)\},\mbox{ or if}\\
&(iii)\,d(a,b) = d(c,d)\mbox{ and }\min\{\mathcal{R}(a),\mathcal{R}(b)\} = \min\{\mathcal{R}(c),\mathcal{R}(d)\}\mbox{ and }\max\{\mathcal{R}(a),\mathcal{R}(b)\} < \max\{\mathcal{R}(c),\mathcal{R}(d)\}.
\end{align*}
\end{definition}

The MST that is constructed by applying Kruskal's algorithm to the edges that are ordered with respect to weight and vertex rank is called a vertex-ranked MST (VRMST).

Now, we will prove Lemma \ref{lem:RelatingMSTToPhylogeneticTree}, which is used to prove the correctness of CLGrouping. 

\begin{lemma}
Adapted from parts $(i)$ and $(ii)$ of Lemma 8 in \cite{Choi2010b}.
Given a phylogenetic tree $T$ and a ranking $\mathcal{R}$ over the labeled vertices in $T$, let $G$ be the distance graph that corresponds to $T=(V_{T},E_{T})$ and let $E_{\leq}$ be the list of edges of $G$ sorted with respect to edge weight and vertex rank, as defined in Definition \ref{def:orderedEdges}. Let $M=(V_{M},E_{M})$ be the VRMST that is constructed by applying Kruskal's algorithm to $E_{\leq}$. The surrogate vertex of each hidden vertex is defined with respect to distance and vertex rank as given in Definition \ref{def:vertexRankBasedSurrogateVertex}. $M$ is related to $T$ as follows.
\label{lem:RelatingMSTToPhylogeneticTree}
\begin{enumerate}[(i)]
\item If $j \in V_{M}$ and $h \in \mbox{Sg}^{-1}(j)$ s.t. $h\neq j$, then every vertex in the path in $T$ that connects $j$ and $h$ belongs to the inverse surrogate set $\mbox{Sg}^{-1}(j)$.
\item For any two vertices that are adjacent in $T$, their surrogate vertices, if distinct, are adjacent in $M$, i.e., for all $i,j\in V_{T}$ with $\mbox{Sg}(i)\neq \mbox{Sg}(j)$, 
$$\{i,j\} \in E_{T} \Rightarrow \{\mbox{Sg}(i),\mbox{Sg}(j)\} \in E_{M}$$
\end{enumerate}
\end{lemma}

\begin{figure}[t]
\centering
\includegraphics[width=0.8\textwidth]{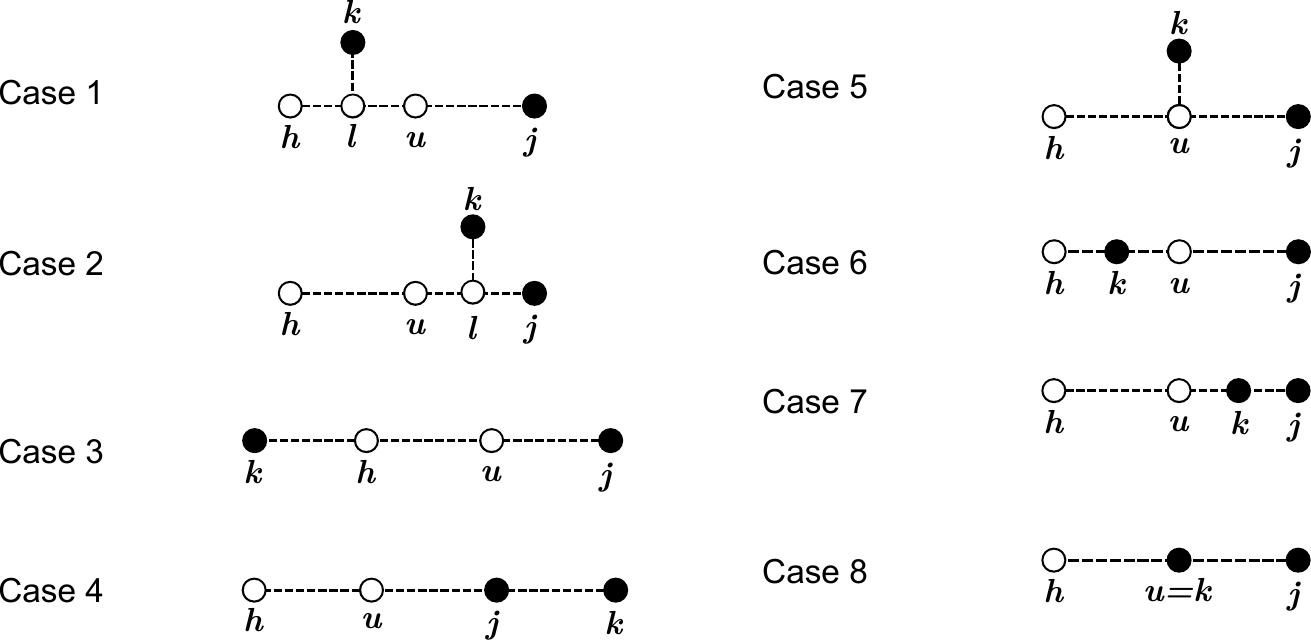}
\caption{The cases that were considered in the proof of Lemma \ref{lem:RelatingMSTToPhylogeneticTree} part $(ii)$. For some phylogenetic tree $T$ let $j$ be a labeled vertex and let $h$ be a hidden vertex in the inverse surrogate set of $j$. $u$ is a vertex in the path between $h$ and $j$. Each case specifies one of the eight possible positions of a labeled vertex $k$ w.r.t $h,u$, and $j$. Hidden vertices are represented with white circles and labeled vertices are represented with black circles. Each dashed line represents a path between the two vertices at its end points.}
\label{fig:casesToConsiderInProofOfLemma1Part2}
\end{figure}

\begin{proof}[Proof \ref{prf:RelatingMSTToPhylogeneticTree}]
\refstepcounter{prf}
\label{prf:RelatingMSTToPhylogeneticTree}
First we will prove Lemma \ref{lem:RelatingMSTToPhylogeneticTree} part $(i)$ by contradiction. 

Assume that there is a vertex $u$ on the path between $h$ and $j$, such that Sg$(u) = k \neq j$. We have $d_{uk} \leq d_{uj}$ (equality holds only if $\mathcal{R}(k) < \mathcal{R}(j)$). Similarly, since Sg$(h) = j$, we have $d_{hj} \leq d_{hk}$ (equality holds only if $\mathcal{R}(j) < \mathcal{R}(k)$)
We consider all eight positions of $k$ w.r.t. $h,u,$ and $j$ (see Fig. \ref{fig:casesToConsiderInProofOfLemma1Part2}).

For case 1 we have
\begin{align*}
&d_{hj} \leq d_{hk}\mbox{ (since Sg$(h) = j$)}\\
\Leftrightarrow &d_{hl}+d_{lu}+d_{uj}\leq d_{hl}+d_{hk}\\
\Leftrightarrow &d_{lu}+d_{uj} \leq d_{lk}\\
\Leftrightarrow &d_{uj} < d_{ul}+d_{lk}\\
\Leftrightarrow &d_{uj} < d_{uk} \mbox{ (contradiction since Sg$(u) = k$).}\\
	\end{align*}

For case 2 we have
\begin{align*}
&d_{hj} \leq d_{hk}\mbox{ (equality holds only if $\mathcal{R}(j) < \mathcal{R}(k)$)}\\
\Leftrightarrow &d_{hu}+d_{ul}+d_{lj}\leq d_{hu}+d_{ul}+d_{lk}\\
\Leftrightarrow &d_{ul}+d_{lj}\leq d_{ul}+d_{lk}\\
\Leftrightarrow &d_{uj}\leq d_{uk} \mbox{ (contradiction since Sg$(u)=k$).}
\end{align*}

For case 3 we have
\begin{align*}
&d_{hj} \leq d_{hk}\\
\Leftrightarrow &d_{hu}+d_{uj}\leq d_{hk}\\
\Leftrightarrow &d_{uj}< d_{hk}+{d_{hu}}\\
\Leftrightarrow &d_{uj}< d_{uk} \mbox{ (contradiction since Sg$(u)=k$).}
\end{align*}

For case 4 we have
\begin{align*}
&d_{uk} = d_{uj} + d_{jk}\mbox{ (see Fig. \ref{fig:casesToConsiderInProofOfLemma1Part2} case 4)}\\
\Leftrightarrow &d_{uk} > d_{uj}\mbox{ (contradiction since Sg$(u)=k$).}
\end{align*}

For case 5 we have
\begin{align*}
&d_{hj} \leq d_{hk}\mbox{ (equality holds only if $\mathcal{R}(j) < \mathcal{R}(k)$)}\\
\Leftrightarrow &d_{hu}+d_{uj} \leq d_{hu} + d_{uk}\\
\Leftrightarrow & d_{uj} \leq d_{uk}\mbox{ (contradiction since Sg$(u)=k$).}
\end{align*}

For cases 6,7, and 8, we have
\begin{align*}
&d_{hj} \leq d_{hk}\\
\Leftrightarrow &d_{hk}+d_{kj} \leq d_{hk}\\
\Leftrightarrow &d_{hk}<d_{hk}\mbox{ (contradiction).}
\end{align*}
Now we will prove part $(ii)$ of Lemma \ref{lem:RelatingMSTToPhylogeneticTree}. Consider the edge $\{i,j\}$ in $E_{T}$ such that $\mbox{Sg}(i)\neq \mbox{Sg}(j)$. Let $V_{i}$ and $V_{j}$ be the sides of the split that is induced by the edge $\{i,j\}$, such that $V_{i}$ and $V_{j}$ contain $i$ and $j$, respectively. Let $L_{i}$ and $L_{j}$ be sets of labeled vertices that are defined as $V_{i}\cap V_{M}$ and $V_{j}\cap V_{M}$ respectively. From part $(i)$ of Lemma \ref{lem:RelatingMSTToPhylogeneticTree} we know that Sg$(i)\in L_{i}$ and Sg$(j)\in L_{j}$. Additionally, for any $k\in L_{i}\backslash\{\mbox{Sg}(i)\}$ and $l\in L_{j}\backslash\{\mbox{Sg}(j)\}$, from the definition of surrogate vertex it follows that
\begin{align*}
d_{ki} &\geq d_{\mbox{Sg}(i)i} (\mbox{equality holds only if } \mathcal{R}(\mbox{Sg}(i))<\mathcal{R}(k))\\
d_{lj} &\geq d_{\mbox{Sg}(j)j} (\mbox{equality holds only if } \mathcal{R}(\mbox{Sg}(j))<\mathcal{R}(l))\\
d_{kj} &= d_{ki}+d_{ij}+d_{lj}\\
&\geq d_{\mbox{Sg}(i)i}+d_{ij}+d_{\mbox{Sg}(j)j}\\
&=d_{\mbox{Sg}(i)\mbox{Sg}(j)}.
\end{align*}

It is clear that 
\begin{equation}
\label{eqn:rankIneq}
\min\{\mathcal{R}(k),\mathcal{R}(l)\} > \min\{\mathcal{R}(\mbox{Sg}(i)),\mathcal{R}(\mbox{Sg}(j))\},
\end{equation}
and that 

\begin{equation}
\label{eqn:surrogateIneq}
d_{kl} \geq d_{\mbox{Sg}(i)\mbox{Sg}(j)}.
\end{equation}

The cut property of MSTs states that given a graph $G=(V,E)$ for each pair $V_{1},V_{2}$ of disjoint sets such that $V_{1}\cup V_{2} = V$, each MST of $G$ contains one of the smallest edges (w.r.t. edge weight) which have one end-point in $V_{1}$ and the other end-point in $V_{2}$. 

Note that the vertex-ranked MST $M$ is constructed using edges that are sorted w.r.t. edge weight and the vertex rank $\mathcal{R}$. From equations (\ref{eqn:rankIneq}) and (\ref{eqn:surrogateIneq}) it is clear that among all edges with one end point in $L_{i}$ and the other end-point in $L_{j}$, the edge $\{\mbox{Sg}(i),\mbox{Sg}(j)\}$ is the smallest edge w.r.t edge weight and vertex rank (see Definition \ref{def:orderedEdges}). Since $L_{i}$ and $L_{j}$ are disjoint sets and $L_{i}\cup L_{j}=V_{M}$, it follows that $\{\mbox{Sg}(i),\mbox{Sg}(j)\} \in E_{M}$.

\end{proof}

CLGrouping can be shown to be correct using Lemma \ref{lem:RelatingMSTToPhylogeneticTree} and the rest of the proof that was provided by \cite{Choi2010b}. Thus if the distances are additive in the model tree, CLGrouping will provably reconstruct the model tree provided that the MST that is used by CLGrouping is a vertex-ranked MST (VRMST).

The authors of CLGrouping provide a matlab implementation of their algorithm. Their implementation reconstructs the model tree even if there are multiple MSTs in the underlying distance graph. The authors' implementation takes as input a distance matrix which has the following property: the row index, and the column index of each labeled vertex is equal. The MST that is constructed in the authors implementation is a vertex-ranked MST, with the rank of each vertex being equal to the corresponding row index of the labeled vertex. We implemented their algorithm in python with no particular order over the input distances and were surprised to find out that the reconstructed tree differed from the model tree, even if the input distances were additive in the model tree.

Depending on the phylogenetic tree, there may be multiple corresponding vertex-ranked MSTs with vastly different numbers of leaves. In the next section we discuss the impact of the number of leaves in a vertex-ranked MST, on the efficiency of parallel implementations of CLGrouping.

\section{Relating the number of leaves in a VRMST to the optimality of the VRMST in the context of CLGrouping}\label{sec:numberOfLeavesInVRMST}
\begin{figure}[h]
\centering
\includegraphics[width=\textwidth]{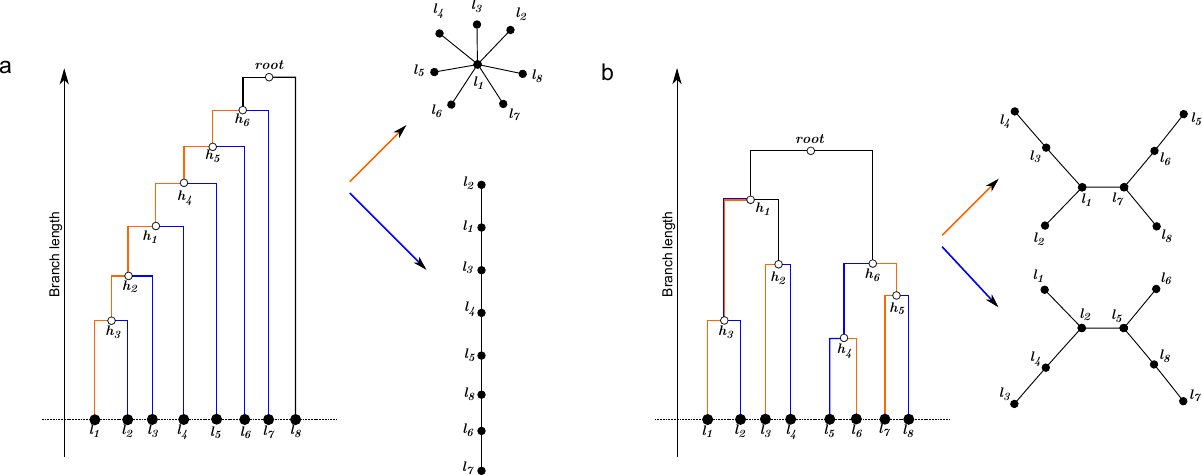}
\caption{Both panels show clock-like phylogenetic trees and VRMSTs with the maximum and the minimum number of leaves, that are constructed by contracting corresponding edges that are highlighted in orange and blue, respectively. The difference between the maximum and the minimum number of leaves in VRMSTs is largest for the caterpillar tree shown in panel a, and smallest for the maximally balanced tree shown in panel b.}
\label{fig:clockLikeTreesAndVRMSTs}
\end{figure}

In the context of parallel programming, \cite{Huang2014} showed that it is possible to parallelize CLGrouping by independently constructing phylogenetic trees for each vertex group, and later combining them in order to construct the full phylogenetic tree.

In order to relate the balancedness of a phylogenetic tree to the number of leaves in a corresponding vertex-ranked MST, we consider clock-like caterpillar trees and maximally balanced trees such that each hidden vertex of each tree has degree three.

Consider the case in which the phylogenetic tree is a caterpillar tree (least balanced). There exists a corresponding VRMST which has a star topology that can be constructed by contracting edges between each hidden vertex and one labeled vertex that is in the surrogate vertex set of each hidden vertex (see Fig \ref{fig:clockLikeTreesAndVRMSTs}a). A star-shaped VRMST has only one vertex group, comprising all the vertices in the VRMST, and does not afford any parallelism.

Instead, if the VRMST was to be constructed by contracting edges between each hidden vertex $h$ and a labeled vertex that is incident to $h$, then the number of the vertex groups would be $n-2$, where $n$ is the number of vertices in the phylogenetic tree. The resulting VRMST would have the minimum number of leaves (two).

With respect to parallelism, an optimal vertex-ranked MST for CLGrouping is a vertex-ranked MST with the maximum number of vertex groups, and equivalently, the minimum number of leaves.

Consider a phylogenetic tree $T=(V_{T},E_{T})$ which is maximally balanced. It is clear that the set $\mathcal{L}(T)$ of labeled vertices of $T$ can be partitioned into a disjoint set $\mathcal{C}$ of vertex pairs such that for each vertex pair $\{u,v\} \in \mathcal{C}$, $u$ and $v$ are adjacent to the same hidden vertex $h\in V_T$. Given a vertex ranking $\mathcal{R}$, the surrogate vertex of $h$ will be $\max_{l \in \{u,v\}} \mathcal{R}(l)$. Thus, independently of vertex ranking, the number of distinct surrogate vertices will be $\mathcal{L}(T)/2$. Each labeled vertex that is not selected as a surrogate vertex will be a leaf in the vertex-ranked MST. It follows that all corresponding VRMSTs of $T$ will have $\mathcal{L}(T)/2$ leaves (see Fig \ref{fig:clockLikeTreesAndVRMSTs}b).

Whether or not the phylogenetic trees that are estimated from real data are clock-like depends on the set of taxa that are being studied. Genetic sequences that are sampled from closely related taxa have been estimated to undergo substitutions at a similar rate, resulting in clock-like phylogenetic trees \citep{DosReis2016}. In the context of evolution, trees are caterpillar-like if there is a strong selection; the longest path from the root represents the best-fit lineage.

In the next section we will present an algorithm for constructing a vertex-ranked MST with the minimum number of leaves. 
\section{Constructing a vertex-ranked MST with the minimum number of leaves}
We aim to construct a vertex-ranked MST with the minimum number of leaves (MLVRMST) from a distance graph. An algorithm for constructing a MLVRMST is presented in subsection \ref{ssec:Algo}. In the following two subsections we will present two lemmas, which will be used for proving the correctness of the algorithm.

\subsection{A common structure that is shared by all MSTs}
In this section we will prove the existence of a laminar family $\mathcal{F}$ over the vertex set of an edge-weighted graph $G$. A collection $\mathcal{F}$ of subsets of a set $S$ is a laminar family over $S$ if, for any two intersecting sets in $\mathcal{F}$, one set contains the other. That is to say, for each pair $S_{1}, S_{2}$ in $\mathcal{F}$ such that $|S_{1}| \leq |S_{2}|$, either $S_{1}\cap S_{2}=\varnothing$, or $S_{1}\subset S_{2}$.

The vertex sets in $\mathcal{F}$ define a structure that is common to each MST of $G$. Furthermore, $\mathcal{F}$ can be used to obtain an upper bound on the degree of each vertex in a MST. The notion of a laminar family has been utilized previously by \cite{Ravi2006}, for designing an approximation algorithm for the minimum-degree MST 

\begin{lemma}
\label{lem:laminarFamily}
Given an edge-weighted graph $G=(V,E)$ with $k$ distinct weight classes $W=\{w_{1},w_{2},\ldots,w_{k}\}$, and an MST $M$ of $G$, let $F_{i}$ be the forest that is formed by removing all edges in $G$ that are heavier than $w_{i}$. Let $\mathcal{C}_{i}$ be the collection comprising the vertex set of each component of $F_{i}$. Consider the collection $\mathcal{F}$ which is constructed as follows: $\mathcal{F}=\left\{\cup_{i=1}^{k}\mathcal{C}_{i}\right\}\cup V$. The following is true:
\begin{enumerate}[(i)]
\item $\mathcal{F}$ is a laminar family over $V$
\item Each vertex set in $\mathcal{F}$ induces a connected subgraph in each MST of $G$
\end{enumerate}
\end{lemma}

\begin{proof} $(i)$.
Consider any two vertex sets $S_{1}$ and $S_{2}$ in $\mathcal{F}$. Let $w_{1}$ and $w_{2}$ be the weights of the heaviest edges in the subgraphs of $M$ that are induced by $S_{1}$ and $S_{2}$, respectively. Let $F_{1}$ and $F_{2}$ be the forests that are formed by removing all edges in $M$ that are heavier than $w_{1}$ and $w_{2}$, respectively. Let $\mathcal{C}_{1}$ and $\mathcal{C}_{2}$ be the collections comprising the vertex set of each component in $F_{1}$ and $F_{2}$, respectively. 

It is clear that $S_{1} \in \mathcal{C}_{1}$ and $S_{2} \in \mathcal{C}_{2}$. Consider the case where $w_{1}=w_{2}$. Since $\mathcal{C}_{1}$=$\mathcal{C}_{2}$, it follows that $S_{1}\cap S_{2}=\varnothing$. If $w_{1}\neq w_{2}$, then without loss of generality, let $w_1 < w_2$. $F_{2}$ can be constructed by adding to $F_{1}$ all edges in $M$ that are no heavier than $w_{2}$. Each component in $F_{1}$ that is not in $F_{2}$ induces a connected subgraph in exactly one component of $F_{2}$. If $S_{1} \in \mathcal{C}_{1}\cap \mathcal{C}_{2}$ then $S_{1}\cap S_{2}=\varnothing$. Otherwise, if $S_{1} \in \mathcal{C}_{1}\backslash \mathcal{C}_{2}$, then $S_{1}$ is a subset of exactly one set in $\mathcal{C}_{2}$. This implies that either $S_{1}\subset S_{2}$, or $S_{1}\cap S_{2}=\varnothing$. Thus $\mathcal{F}$ is a laminar family over $V$. 

$(ii)$. Let $S_{i}$ be the vertex set of a component in the subgraph $G_{i}$ of $G$ that is created by removing all edges in $G_{i}$ that are heavier than $w_{i}$. It is clear that $S_{i}$ induces a connected subgraph in each minimum spanning forest of $G_{i}$. For each minimum spanning forest there is a corresponding MST of $G$, such that the minimum spanning forest can be constructed by removing from the MST all the edges are heavier than $w_{i}$. It follows that $S_{i}$ induces a connected subgraph in each MST of $G$.
\end{proof}

\subsection{Selecting surrogate vertices on the basis of maximum vertex degree}
\begin{lemma}
We are given a phylogenetic tree $T$, the corresponding distance graph $G=(V,E)$, and the laminar family $\mathcal{F}$ of the distance graph. Let the subgraph $g=(V_{g},E_{g})$ of $G$ contain all edges that are present in at least one MST of $G$. Let $h$ be a hidden vertex in $T$ such that there is a leaf $l$ in $\mathbf{Sg}(h)$, and $h$ is incident to $l$. Let $S_{i}$ be a vertex set in $\mathcal{F}$ and let $w_{i}$ be the corresponding edge weight. Then the following holds:
\label{lem:RelatingSurrogateSetToLaminarFamily}
\begin{enumerate}[(i)]
\item Let $J_{v}$ be the set of all vertices that are incident to vertex $v$ in $g$. Let $S_{v}$ be the smallest sub-collection of $\mathcal{F}$ that covers $J_{v}$ but not $v$. Among all MSTs, the maximum vertex degree $\delta_{\max}(v)$ of $v$ is $|S_{v}|$.
\item $\delta_{\max}(l) \leq \delta_{\max}(v)$ for each vertex $v$ in $\mathbf{Sg}(h)$.
\end{enumerate}
\end{lemma}
\begin{proof}

$(i)$. Let $J_{v}=\{j_{1},j_{2},\ldots,j_{k}\}$ be the set of all vertices that are incident to $v$. Let $M$ be some MST of $G$. Let $\mathcal{S}_{v}=\{S_{1},S_{2},\ldots,S_{m}\}$ be the smallest sub-collection of $\mathcal{F}$ that covers $J_{v}$ and does not include $v$. Let $\mathcal{S}_{v}$ contain a set $S_{i}$ that covers multiple vertices in $J$. Let $j_{1}$ and $j_{2}$ be any two vertices in $S_{i}$. Let $w_{i}$ be the heaviest weight on the path that joins $j_{1}$ and $j_{2}$ in $M$. The edges $\{v,j_{1}\}$ and $\{v,j_{1}\}$ are heavier than $w_{i}$. If they were not, then we would have $v \in S_{i}$. Since $v$, $j_{1}$ and $j_{2}$ are on a common cycle, each MST of $G$ can only contain one of the two edges $\{v,j_{1}\}$, and $\{v,j_{2}\}$. It follows that for each set $S_{i} \in \mathcal{S}_{v}$, each MST can contain at most one edge which is incident to $v$ and to a vertex in $S_{i}$. Thus the maximum number of edges that can be incident to $v$ in any MST is the number of vertex sets in $\mathcal{S}_{v}$, i.e., $\delta_{\max}(v)=|\mathcal{S}_{v}|$.

$(ii)$. Let $J_{l}$ and $J_{v}$ be the set of all vertices that are incident to $l$ and $v$ in $g$, respectively. Let $j\in J_{l}\backslash\mathbf{Sg}(h)$. The weight of the edge $\{j,l\} \in E_{g}$ is given by $d_{jl}$. $d_{jh}>d_{vh}$ since $j\notin\mathbf{Sg}(h)$. Thus $d_{lj}>d_{lv}$, and consequently $v\in J_{l}$. We have $d_{jl} = d_{jh}+d_{hl} = d_{jh}+d_{hv} = d_{jv}$. Consider the MST $M=(V_{M},E_{M})$ that contains the edges $\{l,v\}$ and $\{l,h\}$. Consider the spanning tree $M^{\prime}$ that is formed by removing $\{l,h\}$ from $E_{M}$ and adding $\{v,h\}$. $M^{\prime}$ and $M$ have the same sum of edge weights. Thus we also have $j\in J_{v}$. Consequently $J_{l} \subseteq J_{v}$. Let $\mathcal{S}_{l}$ and $\mathcal{S}_{v}$ be the smallest sub-collections of $\mathcal{F}$ such that $\mathcal{S}_{l}$ covers $J_{l}$ but does not contain $l$, and $\mathcal{S}_{v}$ covers $J_{v}$ but does not contain $v$. $\mathcal{S}_{v}$ covers both $J_{l}$ and $J_{v}$ since $J_{l}\subseteq J_{v}$. Thus $|\mathcal{S}_{l}|\leq|\mathcal{S}_{v}|$. From part $(i)$, we know that $|\mathcal{S}_{l}|=\delta_{\max}(l)$ and $|\mathcal{S}_{v}|=\delta_{\max}(v)$. Thus $\delta_{\max}(l)\leq \delta_{\max}(v)$.
\end{proof}

\newpage
\subsection{Constructing a minimum leaves vertex-ranked MST}\label{ssec:Algo}
We now give an overview of Algo. \ref{algo:MLVRMSTFromDistanceGraph}.  Algo. \ref{algo:MLVRMSTFromDistanceGraph} takes as input a distance graph $G = (V,E)$ and computes $\delta_{\max}$ for each vertex in $V$. Subsequently, a ranking $\mathcal{R}$ over $V$ is identified such that vertices with lower $\delta_{\max}$ are assigned higher ranks. The output of Algo. \ref{algo:MLVRMSTFromDistanceGraph} is the vertex-ranked MST which is constructed using $\mathcal{R}$. If $G$ is weighted with tree-additive distances then the output of Algo. \ref{algo:MLVRMSTFromDistanceGraph} is a vertex-ranked MST with the minimum number of leaves (MLVRMST). 

\begin{figure}[htbp]
\centering
\includegraphics[width=0.6\textwidth]{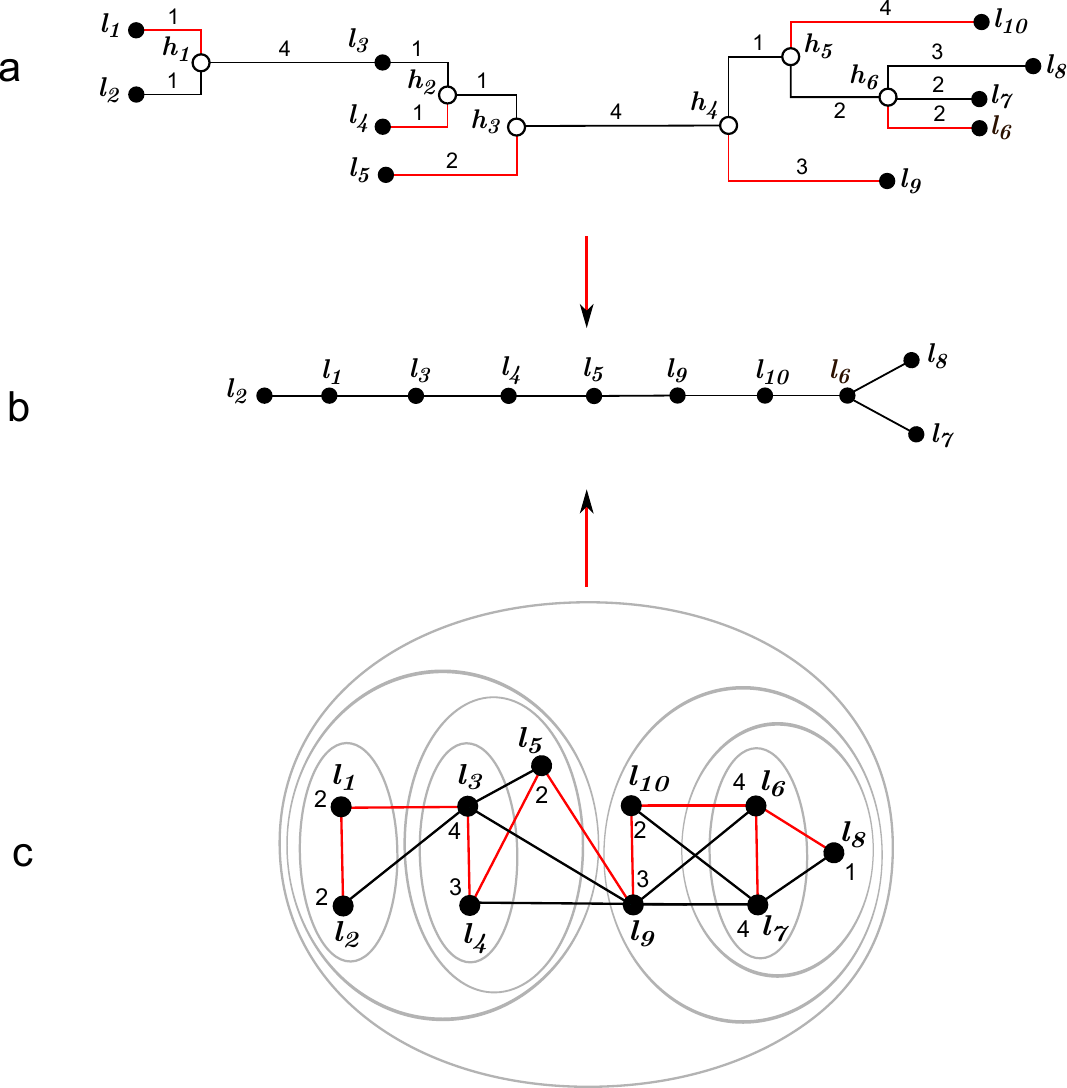}
\caption{Panel a shows a generally labeled phylogenetic tree $T$. Algo. \ref{algo:MLVRMSTFromDistanceGraph} was applied to the distance graph $G$ of $T$. Panel b show the output $M$ of Algo. \ref{algo:MLVRMSTFromDistanceGraph} which is a MLVRMST of $G$. Panel c shows $M$ (in red) superimposed with the laminar family $\mathcal{F}$, and the graph $g$ which contains all edges of $G$ that are present in at least one MST. Additionally each vertex has been labeled with the corresponding $\delta_{\max}$.}
\label{fig:unionGraph}
\end{figure}

An example of a phylogenetic tree, a corresponding MLVRMST, and the output MST $M$ of Algo. \ref{algo:MLVRMSTFromDistanceGraph}, is shown in Fig. \ref{fig:unionGraph}. $M$ is superimposed with the following: the laminar family $\mathcal{F}$, the subgraph $g$, and $\delta_{\max}$ for each vertex.

\begin{algorithm}[H]
\label{algo:MLVRMSTFromDistanceGraph}
\SetKwInput{Initialize}{Initialize}
\KwIn{$G=(V,E)$}
\Initialize{\\
For each vertex in $V$, create a Make-Set object\;
Create empty arrays called $E_{w},W,E_{\textit{fixed}}, E_{\textit{flexible}}$, and $E_{\textit{selected}}$\;
Create empty hash tables called $\textit{CompNbrs}$ and $\textit{CompGraphs}$\;
Create a hash table called $\delta_{\max}$ and for each $u$ in $V$ set $\delta_{\max}(u)$ to zero\;
}
$E_{\leq} \leftarrow$ array of edges of $G$ that are sorted in order of increasing weight\;
$w_{\textit{old}} \leftarrow$ weight of the lightest edge\;
\For{$\{u,v\}$ in $E_{\leq}$}{ \label{line:unionStart}
$w \leftarrow$ weight of $\{u,v\}$\; 
\If{$w > w_{\textit{old}}$}{
$V_{\textit{flexible}} \leftarrow Set()$\;\label{line:start}
Add $w$ to $W$\;
\For{$\{u,v\}$ in $E_{w}$}{
\uIf{\textit{Find}$(u)\neq$ \textit{Find}$(v)$}{
Union($u$,$v$)\;
}
\Else{
Add $u$ to the set $V_{\textit{flexible}}$\; \label{line:unionEnd}
}
}
$V_{w} \leftarrow$ vertices in $E_{w}$\;
\For {$u$ in $V_{w}$}{
Increase $\delta_{\max}(u)$ by $|\textit{CompNbrs}(u)|$\;\label{line:computeDeltaMax}
}
$\textit{CompNbrs} \leftarrow$ empty hash table\;
Add the set $\{u,v,\textit{Comp}(u,w),\textit{Comp}(v,w)\}$ to the array $\textit{CompGraphs}(\textit{Find}(u))$, for each edge $\{u,v\}$ in $E_{w}$ such that $\textit{Find}(u$) equals $\textit{Find}(x$) for at least one vertex $x$ in $V_{\textit{flexible}}$\;
Add to $E_{\textit{fixed}}$, all the edges in $E_{w}$ that are not in $E_{\textit{flexible}}$\;\label{line:end}
$E_{w} \leftarrow$ empty array\;
$w_{\textit{old}}\leftarrow w$\;
}
\If{$Find(u$)$\neq$ $Find(v$)}{\label{line:addEdgesToEw}
Add $\{u,v\}$ to $E_{w}$\;
Set $\textit{Comp}(u,w)$ and \textit{Comp}$(v,w)$ to \textit{Find}($u$) and \textit{Find}($v$), respectively\;
Add $\textit{Comp}(u,w)$ to the set $\textit{CompNbr}(v$)\;
Add $\textit{Comp}(v,w)$ to the set $\textit{CompNbr}(u$)\;
}
}
If $|E_{w}|$ is greater than zero, then repeat lines \ref{line:start} through \ref{line:end}\;
Identify a ranking $\mathcal{R}$ over $V$ such that vertices with lower $\delta_{\max}$ are assigned higher ranks\;
\For{$\textit{CompName}$ in $\textit{CompGraphs}$.keys()}{
$E_{\textit{comp}} \leftarrow$ all edges $\{\textit{Comp}(u,w),\textit{Comp}(v,w)\}$ in $\textit{CompGraphs}(\textit{CompName})$\;
To each vertex \textit{Comp}$(u,w)$ in $E_{\textit{comp}}$, assign the rank $\mathcal{R}(u)$\;
$E_{\leq \textit{comp}} \leftarrow$ edges in $E_{\textit{comp}}$ sorted w.r.t. vertex rank (see Definition \ref{def:orderedEdges}; each edge has weight $w$)\; 
Add to $E_{\textit{selected}}$, each edge in the graph that is constructed by applying Kruskal's algorithm to $E_{\leq \textit{comp}}$\;\label{line:Kruskals}
}
\KwOut{$M=(V,E_{\textit{fixed}}\cup E_{\textit{selected}})$}
\caption{MinLeavesVertexRankedMST of $G$}
\end{algorithm}
\vspace{1 em}
First we prove the correctness of Algo. \ref{algo:MLVRMSTFromDistanceGraph}, and subsequently, we derive its time complexity. Algo. \ref{algo:MLVRMSTFromDistanceGraph} makes use of the disjoint-set data structure, which includes the operations: Make-Set, Find, and Union. The data structure is stored in memory in the form of a forest with self-loops and directed edges. Each directed edge from a vertex points to the parent of the vertex. A Make-Set operation creates a singleton vertex that points to itself. Each component in the forest has a single vertex that points to itself. This vertex is called the root. A Union operation takes as input, the roots of two components, and points one root to the other. A Find operation takes as input a vertex, and returns the root of the component that contains the vertex. Specifically, we implemented balanced Union, and Find with path compression. For a more detailed description please read the survey by \cite{Galil1991}. 

\begin{thm}
Given as input a distance graph such that the distances are additive in some phylogenetic tree with strictly positive branch lengths, Algo. \ref{algo:MLVRMSTFromDistanceGraph} constructs a vertex-ranked MST with the minimum number of leaves.
\end{thm}

\begin{proof}
Let $T=(V_{T},E_{T})$ be the phylogenetic tree that corresponds to the distance graph $G=(V,E)$. Let $W$ be the set of weights of edges in $E$. Let $\mathcal{F}$ be the laminar family over $V$, as defined in Lemma \ref{lem:RelatingSurrogateSetToLaminarFamily}. Let $g$ be the subgraph of $G$ that contains the edges that are present in at least one MST of $G$. Let $M$ be the output of Algo. \ref{algo:MLVRMSTFromDistanceGraph}.

Each edge in $E_{w}$ is incident to vertices in different components. Since edges in $E$ are visited in order of increasing weight, each edge in $E_{w}$ is present in at least one MST of $G$. 

Let $c$ be the root of the component that is formed after Union operations are performed on each edge in $E_{w}$. Let $E_{c}$ be the subset of $E_{w}$ such that each edge in $E_{c}$ is incident to vertices that are in component $c$ after all Union operations on $E_{w}$ have been performed. Let $\mathcal{C}$ be the set of components such that each vertex in $E_{c}$ is contained in a component in $\mathcal{C}$ before any Union operations on $E_{w}$ have been performed. Define the component graph $G_{\mathcal{C}}$ over $\mathcal{C}$ to be the graph whose vertices are elements in $\mathcal{C}$, and whose edges are given by elements in $E_{c}$. It is clear that $G_{\mathcal{C}}$ is connected. We now consider the time point after all Union operations on $E_{w}$ have been performed.

If $G_{\mathcal{C}}$ is a simple graph with no cycles, i.e., $|\mathcal{C}| = |E_{c}|-1$, then each edge in $E_{c}$ must be present in each MST of $G$. All edges in each simple, acyclic, component graph, are stored in $E_{\textit{fixed}}$. If $G_{\mathcal{C}}$ is not simple, or if it contains cycles, then each edge in $G_{\mathcal{C}}$ is stored in $\textit{CompGraphs}(c)$. Additionally each so-called component label $\{\textit{Comp}(u,w),\textit{Comp}(v,w)\}$ is also stored in $\textit{CompGraphs}(c)$. For each vertex $u \in V$ the component label $\textit{Comp}(u,w)$ is the root of the component that contains $u$ before any union operations have been performed on edges in $E_{w}$. For each component $c$, the component graph $G_{\mathcal{C}}$ is induced by the component edges.

Let $\mathcal{S}$ be the smallest sub-collection of the laminar family $\mathcal{F}$ such that $\mathcal{S}$ covers the neighbors of $u$ but not $u$. Let $F_{w}$ be the subgraph of $G$ that is formed by removing from $G$ all edges that are heavier than $w$. Let $\mathcal{N}_{w}$ be the set of vertices in $E_{w}$ that are adjacent to $u$. Let $\mathcal{C}_{w}$ be the collection comprising the vertex set of each component of $F_{w}$ that contains at least one vertex in $\mathcal{N}_{w}$. It is easy to see that $\mathcal{C}_{w} \subset \mathcal{S}$. It follows that $\mathcal{S}=\cup_{w\in W}\mathcal{C}_{w}$, where $W$ is the set comprising the unique edge weights of $G$. Thus $\delta_{\max}(u) = |\mathcal{S}| = \sum_{w \in W} |\mathcal{C}_{w}|$. Thus the operations in line \ref{line:computeDeltaMax} correctly compute $\delta_{\max}(u)$.

At this time point all the edges of $G$ have been visited. Subsequently, Algo. \ref{algo:MLVRMSTFromDistanceGraph} selects a vertex ranking $\mathcal{R}$ such that vertices with lower $\delta_{\max}$ are given higher ranks.

Let $E_{\textit{flexible}}$ be the set containing the edges $\{u,v\}$ that are stored in $\textit{CompGraphs}$. Let Kruskal's algorithm be applied to the edges in $E_{\textit{fixed}}\cup E_{\textit{flexible}}$ that are sorted with respect to weight and $\mathcal{R}$, and let the resulting MST be the vertex-ranked MST $M_{\mathcal{R}} = (V_{\mathcal{R}},E_{\mathcal{R}})$.

Let $S$ be the set of all vertices in $\textit{Comp}(u,w)$. From Lemma \ref{lem:laminarFamily} $(ii)$, we know that $S$ induces a connected subgraph in each MST of $G$. This implies that, after all the edges that are no heavier than $w$ have been visited by Algo. \ref{algo:MLVRMSTFromDistanceGraph}, the vertex set of the component that contains $u$ is independent of the notion of the vertex rank that is used to sort the edges. Thus, instead of applying Kruskal's algorithm to each edge in $E_{\textit{fixed}}\cup E_{\textit{flexible}}$, we can avoid redundant computations by applying Kruskal's algorithm independently to each component graph. Consequently, $E_{\mathcal{R}}= E_{\textit{fixed}}\cup E_{\textit{selected}}$.

From Lemma \ref{lem:RelatingSurrogateSetToLaminarFamily} $(ii)$, we know that, if there is a leaf $l$ in $\mathbf{Sg}(h)$, such that $\{h,l\}\in E_{T}$, then among all vertices in $\mathbf{Sg}(h)$, $\delta_{\max}(l)$ is smallest. Consequently $l$ has the highest rank in $\mathcal{R}$, when compared to other vertices in $\mathbf{Sg}(h)$. Since the surrogate vertex of $h$ is the highest-ranked vertex in $\mathbf{Sg}(h)$, Algo. \ref{algo:MLVRMSTFromDistanceGraph} implicitly selects $l$ as the surrogate vertex of $h$. Since each leaf in $T$ is adjacent to at most one hidden vertex, the vertex ranking that is selected by Algo. \ref{algo:MLVRMSTFromDistanceGraph}, maximizes the number of distinct leaves that are selected as surrogate vertices. Contracting the path in $T$ between a hidden vertex and the corresponding surrogate vertex, increases the degree of the surrogate vertex. Thus, among all vertex-ranked MSTs, $M$ has the minimum number of leaves.
\end{proof}

\subsection{Time complexity of Algorithm
\ref{algo:MLVRMSTFromDistanceGraph}}\label{ss:implementationDetails}
We partition the operations of Algo. \ref{algo:MLVRMSTFromDistanceGraph} into three parts. Part $(i)$ sorts all the edges in $E$ and performs Find and Union operations in order to select the edges in $E_{\textit{fixed}}$ and $\textit{CompGraphs}$. Part $(ii)$ computes $\delta_{\max}$ for each vertex in $V$, and part $(iii)$ sorts, and applies Kruskal's algorithm to the edges in each component graph in $\textit{CompGraphs}$.

In part $(i)$ Algo. \ref{algo:MLVRMSTFromDistanceGraph} iterates over the edges in $G$ which are sorted w.r.t. edge weight. $G=(V,E)$ is a fully connected graph with $n$ vertices and $n(n-1)/2$ edges. We used python's implementation of the Timsort algorithm \citep{Peters2002} which sorts the edges in $O(n^{2}\log n)$ time. Let $m_{f}$ be the number of edges in $E_{\textit{fixed}}$, and let $m_{c}$ be the number of edges that are in a component graph. It is clear that $m_{f}+m_{c} \leq n(n-1)/2$. Algo. \ref{algo:MLVRMSTFromDistanceGraph} iterates over each edge in $G$ and performs $n(n-1)/2$ + $m_{f}$ + $m_{c}$ Find operations, and $n-1$ Union operations. Since we implemented balanced Union, and Find with path compression, the time-complexity of these operations is $O((n(n-1)/2+m_{f}+m_{c})(\alpha((n(n-1)/2+m_{f}+m_{c},n))$ = $O(n^{2}(\alpha((n(n-1)/2+m_{f}+m_{c},n))$, where $\alpha((n(n-1)/2+m_{f}+m_{c},n)$ is the inverse of Ackermann's function as defined in \cite{Tarjan1975}, and is less than 5 for all practical purposes. The total time complexity of part $(i)$ is $O(n^{2}\log n)$.

The operations in line \ref{line:computeDeltaMax} compute $\delta_{\max}(u)$ by counting the number of distinct components that cover the vertices $J_{u}\subset E_{w}$, such that each vertex $j\in J_{u}$ is adjacent to $u$. Assuming that the insertion and retrieval operations on hash tables, and insertion operations arrays have linear time-complexity, the total time complexity of part $(ii)$ is $O(m_{f}+m_{c})$.

Let the number of component graphs in $\textit{CompGraphs}$ be $k$ and let the number of edges and vertices in the $i^{th}$ component graph be $m_{i}$ and $n_{i}$, respectively. The time complexity of sorting, and applying Kruskal's algorithm to $m_{i}$ edges, is $O(m_{i}\log m_{i}) + O(m_{i}\alpha(m_{i},n_{i})) = O(m_{i}\log m_{i})$. The total time complexity of part $(iii)$ is 

\begin{align*}
&\sum_{i=1}^{k}O(m_{i}\log m_{i})\\
&=O\left(\sum_{i=1}^{k}m_{i}\log m_{i}\right)\\
&=O\left(\left(\sum_{i=1}^{k}m_{i}\right)\sum_{i=1}^{k}\dfrac{m_{i}}{\left(\sum_{i=1}^{k}m_{i}\right)}\log m_{i}\right)\\
&=O\left(m_{c}\sum_{i=1}^{k}\dfrac{m_{i}}{m_{c}}\log m_{i}\right)\\
&\leq O\left(m_{c}\log \sum_{i=1}^{k}\dfrac{m_{i}^{2}}{m_{c}}\right)
\quad \text{from Jensen's inequality}\\
&\leq O\left(m_{c}\log \sum_{i=1}^{k}m_{i}^{2}\right)\\
&\leq O\left(m_{c}\log\left(\sum_{i=1}^{k}m_{i}\right)^{2}\right)\\
&=O(m_{c}\log m_{c})
\end{align*}

The total time complexity of Algo. \ref{algo:MLVRMSTFromDistanceGraph} is $O(n^{2}\log n) + O(m_{f}+m_{c}) + O(m_{c}\log m_{c}) = O(n^{2}\log n)$.

\section{Computational complexity of the MLVRMST construction problem}
Let $\mathcal{T}$ be the set of all phylogenetic trees. Let $\mathcal{G}$ be the set of edge-weighted graphs, such that the edges of each graph in $\mathcal{G}$ are weighted with distances that additive in some tree in $\mathcal{T}$. Algo. \ref{algo:MLVRMSTFromDistanceGraph} constructs a MLVRMST of any graph in $\mathcal{G}$, in time $O(n^{2}\log n)$. Thus, for graphs in $\mathcal{G}$, the decision version of the optimization problem MLVRMST is in the complexity class \textbf{P}. For graphs whose edges are not weighted with tree-additive distances, the MLVRMST problem may not be in \textbf{P}.

Consider the general optimization problem of constructing an MST with the minimum number of leaves (MLMST). Since the decision version of MLMST can be verified in polynomial time, MLMST is in \textbf{NP}. Additionally, it is easy to show that there is a polynomial time reduction from the Hamiltonian path problem to MLMST. Since the Hamiltonian path problem is in \textbf{NP-complete}, MLMST must be in $\textbf{NP-hard} \cap \textbf{NP}=\textbf{NP-complete}$.

\section{Acknowledgements}
We thank Erik Jan van Leeuwen and Davis Isaac for helpful discussions during the early stages
of the work presented here.

\section{Funding}
PK's work has been funded in part by the German Center for Infection Research (DZIF, German Ministry of Education and Research Grants No. TTU 05.805, TTU 05.809).

\section{Availability of code}
A python implementation of Algo. \ref{algo:MLVRMSTFromDistanceGraph} can be found at\\
http://resources.mpi-inf.mpg.de/departments/d3/publications/prabhavk/minLeavesVertexRankedMST

\bibliographystyle{natbib}
\bibliography{ML-VR-MST}
\end{document}